\documentclass[letterpaper, 10 pt, conference]{ieeeconf}  

\IEEEoverridecommandlockouts                              

\usepackage{graphics} 
\usepackage{hyperref}%
\hypersetup{colorlinks=true, linkcolor=blue, breaklinks=true, urlcolor=blue, citecolor=blue}

\usepackage{amsmath} 
\usepackage{amssymb}  

\usepackage[prependcaption,colorinlistoftodos]{todonotes}

\newcommand{\R}{\mathbb{R}}
\newcommand{\s}{\mathbb{S}}
\usepackage{mathtools}

\newcommand{\map}[3]{#1 \colon #2 \rightarrow #3}

\usepackage{xcolor}

\usepackage{comment}
\usepackage{algorithmic}
\usepackage{algorithm}

\newtheorem{definition}{Definition}
\newtheorem{assumption}{Assumption}
\newtheorem{lemma}{Lemma}

\newtheorem{theorem}{Theorem}

\newtheorem{remark}{Remark}

\newcommand\setItemnumber[1]{\setcounter{enumi}{\numexpr#1-1\relax}}

\allowdisplaybreaks

\title{\LARGE \bf Contractivity Analysis and Control Design for Lur'e
  Systems: \\Lipschitz, Incrementally Sector Bounded, and Monotone
  Nonlinearities}

\author{Ryotaro Shima$^{1}$, Alexander Davydov$^{1}$ and Francesco Bullo$^{1}$ 
\thanks{This work was supported in part by AFOSR award FA9550-22-1-0059.}
\thanks{$^{1}$Ryotaro Shima, Alexander Davydov, and Francesco Bullo are with the Center for Control, Dynamical Systems, and Computation, University of California at Santa Barbara, Santa Barbara, CA 93106 USA.
    {\tt\small \{rshima,davydov,bullo\}@ucsb.edu}}%
}

\begin{document}

\maketitle
\thispagestyle{empty}
\pagestyle{empty}

\begin{abstract}
  In this paper, we study the contractivity of Lur'e dynamical systems
  whose nonlinearity is either Lipschitz, incrementally sector bounded, or
  monotone.  We consider both the discrete- and continuous-time
  settings. In each case, we provide state-independent linear matrix
  inequalities (LMIs) which are necessary and sufficient for contractivity.
  Additionally, we provide LMIs for the design of controller gains such
  that the closed-loop system is contracting.  Finally, we provide a
  numerical example for control design.
\end{abstract}


\section{INTRODUCTION}

Contraction theory has attracted substantial interest within the nonlinear
control community due to its capability to ensure exponential convergence
between arbitrary trajectories of nonlinear systems \cite{WL-JJES:98,
  FB:24-CTDS}. Beyond traditional stability analysis, contractivity
provides significant advantages for the investigation of nonlinear
dynamical phenomena, such as synchronization in networked systems
\cite{LS-MA-EDS:10, FZ-HLT-JMAS:14}, design of distributed controllers
\cite{HSS-MR-IRM:19}, and enhanced robustness to noise and disturbances
\cite{AD-VC-AG-GR-FB:23f, SY-CM-HC-FA:13}.

Contractivity is typically verified through an inequality condition related
to the exponential rate at which the distance between trajectories
decreases. For nonlinear systems, however, this inequality is often
state-dependent, making its verification challenging. One common solution
is the use of a polytopic relaxation, which simplifies the contractivity
condition to checking stability at the vertices of a polytope corresponding
to linear dynamics \cite{CF-JK-XJ-SM:18}. This approach has been
successfully applied in contexts such as safety verification
\cite{CF-JK-XJ-SM:18} and tube-based model predictive control (tube-MPC)
\cite{SY-CB-HC-FA:10}. Nonetheless, polytopic relaxation frequently leads
to conservative conditions, especially for nonlinearities that are
Lipschitz, incrementally sector bounded, or monotone, but not explicitly
polytopic.

Another method involves imposing incremental quadratic constraints on the
nonlinearity. The notable early work in \cite{LDA-MC:13} presents
contractivity conditions assuming that the Lur’e nonlinearity satisfies an
incremental inequality constraint. Similar conditions based on monotonicity
or Lipschitz continuity have been independently established in
\cite{MG-VA-ST-DA:23, MF-MM-GJP:20, MR-RW-IRM:20}. Additionally,
\cite{MG-VA-ST-DA:23} and \cite{SY-CM-HC-FA:13} introduce LMIs for
designing controllers that guarantee the contractivity of closed-loop
systems. More recently, \cite{AD-FB:24i} derived necessary and sufficient
conditions for contractivity, assuming that the nonlinearity is cocoercive.

\paragraph*{Contributions}
In this paper, we generalize previous results to obtain unified necessary
and sufficient conditions for verifying contractivity in Lur’e systems with
Lipschitz, incrementally sector-bounded, and monotone
nonlinearities. Extending existing methods from \cite{MG-VA-ST-DA:23,
  SY-CM-HC-FA:13}, we provide linear matrix inequalities (LMIs) for
designing controllers that guarantee closed-loop contractivity.

Our specific contributions are summarized as follows:\\
1) In Theorem~\ref{thm:lmi_lipschitz}, we introduce an LMI involving the
gain matrices and a weighting norm, providing necessary and sufficient
conditions for the contractivity of closed-loop systems with Lipschitz
nonlinearities. Our LMI generalizes the conditions from
\cite{SY-CM-HC-FA:13} by accommodating a broader class of dynamical systems
and, notably, proving necessity, whereas previous results established only
sufficiency. A detailed comparison is presented in
Subsection~\ref{subsubsection:comparison_lipschitz}.\\
2) We extend the result of Theorem~\ref{thm:lmi_lipschitz} to incrementally
sector-bounded and monotone nonlinearities in
Theorem~\ref{thm:lmi_isbmonotone}. Additionally, we demonstrate that
monotonicity emerges as a special case of incremental sector boundedness
under a mild symmetry assumption. While prior works such as
\cite{MG-VA-ST-DA:23} have established only sufficient conditions, our
formulation provides both necessity and sufficiency and an explicit
estimate of the contractivity rate.\\
3) We derive analogous conditions for discrete-time Lur’e systems,
obtaining corresponding LMIs to verify contractivity. This result
facilitates the direct synthesis and practical implementation of digital
controllers.\\
4) We provide a concrete example illustrating our proposed
controller-design method for a Lipschitz Lur’e system, explicitly
demonstrating the applicability of our LMIs to ensure contractivity.

\paragraph*{Notation}

We denote the $n$-dimensional identity matrix by $I_n$, the set of positive
definite matrices of size $n\times n$ by $\s^n_+$, the weighted Euclidean
norm $\sqrt{v^\top P v}$ by $\|v\|_P$, where $P\in \s^n_+$, and a
block-diagonal matrix whose diagonal blocks are $A$ and $B$ by
$\mathrm{diag}[A, B]$.  $\langle A \rangle \coloneq A + A^\top$, where $A$
is a square matrix. $[v; w] \coloneqq \begin{bmatrix} v^\top & w^\top \end{bmatrix}^\top$, where $v$ and $w$ are vectors. For symmetric matrices $A, B$, we say that $A \succeq
B$ if $A - B$ is positive semidefinite.

\section{PRELIMINARIES}

\subsection{Lur'e Dynamics and Its Contraction}

Consider the following continuous-time Lur'e system (hereinafter we call CTLS):
\begin{subequations} \label{eq:CTLS}
\begin{align}
    \dot{x} &= Ax + B_\Psi\Psi(y)+Bu
    ,
    \label{plant}
    \\
    y &= Cx
    ,
    \label{def_y}
    \\
    u &= Kx + K_\Psi \Psi(y)
    ,
    \label{def_controller}
\end{align}
\end{subequations}
where $\dot{x}$ denotes the time derivative of $x$, $x\in \R^{n_x}$, $y\in \R^{n_y}$, $u \in \R^{n_u}$, $\Psi : \R^{n_y} \to \R^{n_\Psi}$, and $A, B, B_\Psi, C, K, K_\Psi$ are matrices of appropriate sizes.
Hereinafter, we refer to $\Psi$ as a Lur'e nonlinearity.
The closed-loop system is in the following continuous-time Lur'e form:
\begin{align}
    \dot{x} &= A_\mathrm{cl}x + B_\mathrm{cl}\Psi(Cx)
    ,
    \label{closed_loop}
    \\
     A_\mathrm{cl} &\coloneqq A + BK
    ,
    \label{def_Acl}
    \\
    B_\mathrm{cl} &\coloneqq B_\Psi + BK_\Psi
    .
    \label{def_Bcl}
\end{align}

\begin{assumption}
\label{ass:C}
$n_x \geq n_y$ and $C \in \R^{n_y\times n_x}$ is row full rank.
\end{assumption}
\smallskip
One can always make $C$ row full rank by redefining $\Psi$.

\begin{definition}[Contractivity of CTLS]
\label{def:contraction}
Consider $P \in \s_+^{n_x}$ and a positive constant $\eta$.
The CTLS \eqref{eq:CTLS} is said to be strongly infinitesimally contracting with respect to norm $\|\cdot\|_P$ and rate $\eta$ if, for any pair of solutions $x^{(1)}(t), x^{(2)}(t)$ of the CTLS \eqref{eq:CTLS}, the following inequality holds for any $t \geq s \geq 0$:
\begin{align}
    \| x^{(1)}(t) - x^{(2)}(t) \|_P \leq  \mathrm{e}^{-\eta (t-s)}\| x^{(1)}(s) - x^{(2)}(s) \|_P
    .
    \label{contraction}
\end{align}
\end{definition}
\smallskip

\begin{remark}
The CTLS \eqref{eq:CTLS} is strongly infinitesimally contracting with respect to norm $\|\cdot\|_P$ and rate $\eta$
if and only if the following inequality holds for all $x^{(1)}, x^{(2)} \in \R^{n_x}$:
\begin{align}
    \Delta x ^\top P( A_\mathrm{cl} \Delta x + B_\mathrm{cl} \Delta_x \Psi)
    + \eta \Delta x^\top P \Delta x
    \leq 0
    ,
    \label{contraction_incremental}
    \\
    \Delta x \coloneqq x^{(1)} - x^{(2)}, \ \ \Delta_x \Psi \coloneqq \Psi(Cx^{(1)}) - \Psi(Cx^{(2)})
    .
    \label{def_deltax_deltapsi}
\end{align}
In addition, suppose $\Psi$ is differentiable. Then, \eqref{contraction_incremental} hold if and only if the following inequality holds for all $x \in \R^{n_x}$:
\begin{align}
    \Bigl \langle P \Bigl (A_\mathrm{cl} + B_\mathrm{cl} \frac{\partial \Psi}{\partial y}(Cx) C \Bigl ) \Bigl \rangle
    + 2 \eta P
    \preceq 0
    .
    \label{contraction_differential}
\end{align}
See \cite[Theorem 29]{AD-SJ-FB:20o}) for these equivalences.
\end{remark}

In parallel, we consider the following discrete-time Lur'e system (hereinafter we call DTLS):
\begin{subequations} \label{eq:DTLS}
\begin{align}
    x^+ &= Ax + B_\Psi\Psi(y)+Bu
    ,
    \label{plant_discrete}
    \\
    y &= Cx
    ,
    \label{def_y_discrete}
    \\
    u &= Kx + K_\Psi \Psi(y)
    ,
    \label{def_controller_discrete}
\end{align}
\end{subequations}
where $x^+$ denotes $x$ at the next time step.
The dynamics can be rewritten in the following discrete-time Lur'e form:
\begin{align}
    x^+ &= A_\mathrm{cl}x + B_\mathrm{cl}\Psi(Cx)
    ,
    \label{closed_loop_discrete}
\end{align}
where $A_\mathrm{cl}$ and $B_\mathrm{cl}$ are defined as in \eqref{def_Acl} and \eqref{def_Bcl}, respectively.

\begin{definition}[Contractivity of DTLS]
\label{def:contraction_discrete}
Consider $P \in \s_+^{n_x}$ and a constant $0 \leq \eta < 1$.
The DTLS \eqref{eq:DTLS} is said to be strongly contracting with respect to norm $\|\cdot\|_P$ and factor $\eta$ if, for any pair of solutions $x^{(1)}(k), x^{(2)}(k)$ of the DTLS \eqref{eq:DTLS}, the following inequality holds for any $k \geq \ell \geq 0$:
\begin{align}
    \| x^{(1)}(k) - x^{(2)}(k) \|_P \leq \eta^{k-\ell} \| x^{(1)}(\ell) - x^{(2)}(\ell) \|_P
    .
\end{align}
\end{definition}
\smallskip

\begin{remark}
The DTLS \eqref{eq:DTLS} is strongly contracting with respect to norm $\|\cdot\|_P$ and factor $\eta$ if and only if
the following inequality holds for all $x^{(1)}, x^{(2)} \in \R^{n_x}$:
\begin{gather}
    \begin{aligned}
    (A_\mathrm{cl} \Delta x + B_\mathrm{cl} \Delta_x \Psi)^\top
    P
    (A_\mathrm{cl} \Delta x + B_\mathrm{cl} \Delta_x \Psi)
    &
    \\
    - \eta^2 \Delta x ^\top P \Delta x
    \leq 0
    ,&
    \end{aligned}
    \label{contraction_incremental_discrete}
    \\
    \Delta x \coloneqq x^{(1)} - x^{(2)}, \ \ \Delta_x \Psi \coloneqq \Psi(Cx^{(1)}) - \Psi(Cx^{(2)})
    .
    \label{def_deltax_deltapsi_discrete}
\end{gather}
In addition, suppose $\Psi$ is differentiable.
As in the CTLS, \eqref{contraction_incremental_discrete} holds if and only if the following inequality holds for all $x \in \R^{n_x}$:
\begin{align}
\begin{aligned}
    \Big (
    A_\mathrm{cl} + B_\mathrm{cl} \frac{\partial \Psi}{\partial y}(Cx)
    \Big)^\top
    P
    \Big (
    A_\mathrm{cl} + B_\mathrm{cl} \frac{\partial \Psi}{\partial y}(Cx)
    \Big)
    &
    \\
    - \eta^2 P
    &
    \preceq 0
    .
    \end{aligned}
    \label{contraction_differential_discrete}
\end{align}
See \cite[Theorem 3.7]{FB:24-CTDS} for these equivalences.
\end{remark}

In this paper, we discuss how to design the gain matrices $K, K_\Psi$ with which the closed-loop system is contracting.
The difficulty in such gain design lies in the fact that \eqref{contraction_incremental} (\eqref{contraction_differential}, \eqref{contraction_incremental_discrete}, or \eqref{contraction_differential_discrete}) must be satisfied throughout the space of $\R^{n_x}$ due to the nonlinearity $\Psi$.
To address this, we impose a nonlinearity condition on $\Psi$.

\subsection{Nonlinearity Conditions}

In this paper, we consider three types of conditions on the Lur'e nonlinearity: Lipschitzness, incremental sector boundedness, and monotonicity.

\begin{definition}[Lipschitzness]
\label{def:quadratic_constraint}
Let $\Theta_y \in \s_+^{n_y}, \Theta_{\Psi} \in \s_+^{n_\Psi}$.
A function $\map{\Psi}{\R^{n_y}}{\R^{n_{\Psi}}}$ is said to be $\rho$-Lipschitz with input norm $\|\cdot\|_{\Theta_y}$ and output norm $\|\cdot\|_{\Theta_\Psi}$ if the following inequality holds for each $y^{(1)}, y^{(2)} \in \R^{n_y}$:
\begin{align}
    \Delta_y \Psi^\top \Theta_\Psi \Delta_y \Psi
    \leq
    \rho^2
    \Delta y^\top \Theta_y \Delta y
    ,
    \label{Lipschitz}
\end{align}
where $\Delta y \coloneqq y^{(1)} - y^{(2)}$, $\Delta_y \Psi \coloneqq \Psi(y^{(1)}) - \Psi(y^{(2)})$.
\end{definition}
\smallskip

\begin{remark}
If $\Psi$ is differentiable, then $\Psi$ is $\rho$-Lipschitz with input norm $\|\cdot\|_{\Theta_y}$ and output norm $\|\cdot\|_{\Theta_\Psi}$ if and only if for all $y \in \R^{n_y}$, the following inequality holds, 
\begin{align}
    \Bigl (
    \frac{\partial \Psi}{\partial y}(y) 
    \Bigl )^\top
    \Theta_\Psi
    \frac{\partial \Psi}{\partial y}(y)
    \preceq
    \rho^2
    \Theta_y
    .
    \label{Lipschitz_differential}
\end{align} 
\end{remark}
\smallskip

\begin{definition}[incremental sector boundedness]
\label{def:quadratic_constraint}
Consider $\Theta \in \s_+^{n_\Psi}$, $\Gamma \in \R^{n_y \times n_\Psi}$, and a function $\Psi : \R^{n_y} \to \R^{n_\Psi}$.
$\Psi$ is said to be incrementally sector bounded with sector bound $[0, \Gamma]$ and weight $\Theta$ if the following inequality holds for all $y^{(1)}, y^{(2)} \in \R^{n_y}$:
\begin{align}
    \Delta_y \Psi^\top \Theta ( \Delta_y \Psi - \Gamma \Delta y)
    \leq 0
    ,
    \label{incremental_sector_boundedness}
\end{align}
where $\Delta y \coloneqq y^{(1)} - y^{(2)}$ and $\Delta_y \Psi \coloneqq \Psi(y^{(1)}) - \Psi(y^{(2)})$.
Similarly, $\Psi$ is said to be differentially sector bounded with sector bound $[0, \Gamma]$ and weight $\Theta$ if the following inequality holds for all $y \in \R^{n_y}$:
\begin{align}
    \Bigl \langle
    \Bigl (
    \frac{\partial \Psi}{\partial y}(y) 
    \Bigl )^\top
    \Theta
    \Bigl (
    \frac{\partial \Psi}{\partial y}(y)
    - \Gamma
    \Bigl )^\top
    \Bigl \rangle
    \preceq
    0
    .
    \label{rho_cocoercive_differential}
\end{align} 
\end{definition}
\smallskip

\begin{remark}
The differential sector boundedness in the above definition is addressed in
\cite{MG-VA-ST-DA:23}.  Furthermore, the incremental sector boundedness
with sector bound $[0, \frac{1}{\beta} I]$, $\beta>0$, and weight $\Theta =
I$ is called cocoerciveness \cite{EKR-SB:16, AD-FB:24i}.
\end{remark}

\begin{definition}[monotonicity]
    Suppose $n_y=n_\Psi$.
    Let $\Gamma \in \s_+^{n_\Psi}$.
    A differentiable function $\Psi : \R^{n_\Psi} \to \R^{n_\Psi}$ is said to be monotone with upper bound $\Gamma$ if the following property holds for all $y\in \R^{n_y}$:
    \begin{align}
        0 \preceq
        \frac{1}{2}
        \left (
        \frac{\partial \Psi}{\partial y}(y) + \Bigl ( \frac{\partial \Psi}{\partial y}(y) \Bigl ) ^\top
        \right )
        \preceq \Gamma
        .
        \label{differentially_sector_bounded}
    \end{align}
\end{definition}
\smallskip

The monotonicity in the above definition is considered in \cite{MG-VA-ST-DA:23}.
We will later prove that monotonicity with an upper bound is a special case of the incremental sector boundedness under the following assumption.
\begin{assumption}
    \label{ass:symmetric}
    The derivative of $\Psi$ is symmetric, i.e.,
    \begin{align}
        \frac{\partial \Psi}{\partial y}(y) = \Bigl ( \frac{\partial \Psi}{\partial y}(y) \Bigl ) ^\top
        \ \ \ \forall y \in \R^{n_y}
        .
        \label{symmetric}
    \end{align}
\end{assumption}
\smallskip
Note that \cite{MG-VA-ST-DA:23} has also assumed \eqref{symmetric} to derive an LMI for controller design.
Many functions, such as an element-wise activation function, a log-sum-exp function, and the gradient of scalar functions, satisfy \eqref{symmetric}.

\section{MAIN RESULTS}

In this section, we provide contraction analysis of the Lur'e systems and LMIs for controller design which guarantees closed-loop contractivity.
In Subsection \ref{subsection:lipschitz}, we present our results for Lipschitz Lur'e nonlinearity.
In Subsection \ref{subsection:isb_monotone}, we present the counterpart for incrementally sector bounded or monotone Lur'e nonlinearity, together with a lemma claiming that the monotonicity is a special case of the incremental sector boundedness under Assumption~\ref{ass:symmetric}.
Our LMI is compared with that of literature in Subsection \ref{subsection:comparison}.

Note that statement 2) in each theorem provides an LMI with respect to $W, Z, K_\Psi$.
Using a feasible solution $(W, Z, K_\Psi)$ of the LMI, we obtain controller gains $K$ and $K_\Psi$ in \eqref{def_controller} (or \eqref{def_controller_discrete}) where $K = ZW^{-1}$.
In addition, our matrix inequalities are necessary and sufficient for contractivity under each Lur'e nonlienarity.

\subsection{Contraction under Lipschitzness of Lur'e nonlinearity}
\label{subsection:lipschitz}

\begin{theorem}[Continuous-time Lipschitz Lur'e models]
\label{thm:lmi_lipschitz}
Consider the CTLS \eqref{eq:CTLS} with Assumption \ref{ass:C}.
Let $A_\mathrm{cl}$ as in \eqref{def_Acl}, $B_\mathrm{cl}$ as in \eqref{def_Bcl}, $\Theta_y \in \s_+^{n_y}$, $\Theta_\Psi \in \s_+^{n_\Psi}$, and $\eta$ be a positive constant.
Suppose $B_\mathrm{cl} \neq 0$.
Then, the following three statements are equivalent:
\begin{enumerate}
\item There exists $P \in \s_+^{n_x}$ that satisfies
\begin{align}
    \begin{bmatrix}
        \langle P A_\mathrm{cl} \rangle + 2\eta P + \rho^2 C^\top \Theta_y C & P B_\mathrm{cl} \\
        B_\mathrm{cl}^\top P & -\Theta_\Psi
    \end{bmatrix}
    \preceq 0.
    \label{contraction_multiplier_lipschitz}
\end{align}
\item There exists $P \in \s_+^{n_x}$ that satisfies
\begin{align}
    \begin{bmatrix}
        \langle AW + BZ \rangle + 2\eta W & 
        B_\mathrm{cl} &
        W C^\top
        \\
        B_\mathrm{cl}^\top &
        - \Theta_\Psi &
        0
        \\
        C W &
        0 &
        - \frac{1}{\rho^2} \Theta_y^{-1}
    \end{bmatrix}
    \preceq 0
    ,
    \label{contraction_lmi_lipschitz}
\end{align}
where $W \coloneqq P^{-1}, Z \coloneqq K P^{-1}, B_\mathrm{cl} = B_\Psi + BK_\Psi$.
\item There exists $P \in \s_+^{n_x}$ such that the CTLS \eqref{eq:CTLS} is strongly infinitesimally contracting with norm $\|\cdot\|_P$ and rate $\eta$ for each nonlinearity $\Psi$ that is $\rho$-Lipschitz with input norm $\|\cdot\|_{\Theta_y}$ and output norm $\|\cdot\|_{\Theta_\Psi}$.
\end{enumerate}
\end{theorem}
\smallskip

The proof of Theorem \ref{thm:lmi_lipschitz} is shown in Appendix \ref{appendix:proof:lmi_lipschitz}.

\begin{theorem}[Discrete-time Lipschitz Lur'e models]
\label{thm:lmi_lipschitz_discrete}
Consider the DTLS \eqref{eq:DTLS} with Assumption \ref{ass:C}.
Let $A_\mathrm{cl}$ as in \eqref{def_Acl}, $B_\mathrm{cl}$ as in \eqref{def_Bcl}, $\Theta_y \in \s_+^{n_y}$, $\Theta_\Psi \in \s_+^{n_\Psi}$, and $\eta$ be a positive constant satisfying $0<\eta<1$.
Then, the following three statements are equivalent:
\begin{enumerate}
\item There exists $P \in \s_+^{n_x}$ that satisfies
\begin{align}
    \hspace{-15pt}
    \begin{bmatrix}
        A_\mathrm{cl}^\top P A_\mathrm{cl} - \eta^2 P + \rho^2 C^\top \Theta_y C & A_\mathrm{cl}^\top P B_\mathrm{cl} \\
        B_\mathrm{cl}^\top P A_\mathrm{cl} & B_\mathrm{cl}^\top P B_\mathrm{cl} - \Theta_\Psi
    \end{bmatrix}
    \preceq 0.
    \label{contraction_multiplier_lipschitz_discrete}
\end{align}
\item There exists $P \in \s_+^{n_x}$ that satisfies
\begin{align}
    \hspace{-15pt}
    \begin{bmatrix}
        - \eta^2 W & 
        0 &
        W C^\top &
        (AW + BZ)^\top
        \\
        0 &
        - \Theta_\Psi &
        0 &
        B_\mathrm{cl}^\top
        \\
        C W &
        0 &
        - \frac{1}{\rho^2} \Theta_y^{-1} &
        0
        \\
        AW + BZ &
        B_\mathrm{cl} &
        0 &
        - W
    \end{bmatrix}
    \preceq 0
    ,
    \label{contraction_lmi_lipschitz_discrete}
\end{align}
where $W \coloneqq P^{-1}, Z \coloneqq K P^{-1}, B_\mathrm{cl} = B_\Psi + BK_\Psi$.
\item There exists $P \in \s_+^{n_x}$ such that the system \eqref{closed_loop_discrete} is strongly contracting with norm $\|\cdot\|_P$ and factor $\eta$ for each nonlinearity $\Psi$ that is $\rho$-Lipschitz with input norm $\|\cdot\|_{\Theta_y}$ and output norm $\|\cdot\|_{\Theta_\Psi}$.
\end{enumerate}
\end{theorem}
\smallskip

The proof of Theorem \ref{thm:lmi_lipschitz_discrete} is shown in Appendix \ref{appendix:proof:lmi_lipschitz_discrete}.

\subsection{Contraction under incremental sector boundedness or monotonicity of Lur'e nonlinearity}
\label{subsection:isb_monotone}

\begin{lemma}
    \label{lemma:equiv_isb_monotonic}
    Consider a function $\Psi : \R^{n_\Psi} \to \R^{n_\Psi}$ satisfying Assumption \ref{ass:symmetric}.
    Let $\Gamma \in \s_+^{n_\Psi}$.
    Then, $\Psi$ is monotone with an upper bound $\Gamma$ if and only if $\Psi$ is differentially sector bounded with sector bound $[0, \Gamma]$ and weight $\Gamma^{-1}$.
\end{lemma}
\smallskip

The proof of Lemma \ref{lemma:equiv_isb_monotonic} is shown in Appendix \ref{appendix:proof:equiv_isb_monotone}.

\begin{theorem}[Continuous-time sector-bounded Lur'e models]
\label{thm:lmi_isbmonotone}
Consider the CTLS \eqref{eq:CTLS} with Assumption \ref{ass:C}.
Let $A_\mathrm{cl}$ as in \eqref{def_Acl}, $B_\mathrm{cl}$ as in \eqref{def_Bcl}, $\Theta \in \s_+^{n_\Psi}$, $\Gamma \in \R^{n_y \times n_\Psi}$, $G \coloneqq C^\top \Gamma^\top \Theta$, and $\eta$ be a positive constant.
Then, the following three statements are equivalent:
\begin{enumerate}
\item There exists $P \in \s_+^{n_x}$ that satisfies
\begin{align}
    \begin{bmatrix}
        \langle P A_\mathrm{cl} \rangle + 2\eta P & P B_\mathrm{cl} + G \\
        B_\mathrm{cl}^\top P + G^\top & - 2 \Theta
    \end{bmatrix}
    \preceq 0.
    \label{contraction_multiplier}
\end{align}
\item There exists $P \in \s_+^{n_x}$ that satisfies
\begin{gather}
    \begin{bmatrix}
        \langle AW + BZ \rangle + 2\eta W & 
        B_\mathrm{cl} + WG
        \\
        (B_\mathrm{cl} + WG)^\top &
        -2 \Theta
    \end{bmatrix}
    \preceq 0
    ,
    \label{contraction_lmi}
\end{gather}
where $W \coloneqq P^{-1}, Z \coloneqq K P^{-1}, B_\mathrm{cl} = B_\Psi + BK_\Psi$.
\item There exists $P \in \s_+^{n_x}$ such that the CTLS \eqref{eq:CTLS} is strongly infinitesimally contracting with norm $\|\cdot\|_P$ and rate $\eta$ for each nonlinearity $\Psi$ that is incrementally sector bounded with sector bound $[0, \Gamma]$ and weight $\Theta$.
\end{enumerate}
Furthermore, if $\Psi$ is differentiable, the above three statements are equivalent to the following statement.
\begin{enumerate}
\setItemnumber{4}
\item There exists $P \in \s_+^{n_x}$ such that the CTLS \eqref{eq:CTLS} is strongly infinitesimally contracting with norm $\|\cdot\|_P$ and rate $\eta$ for each nonlinearity $\Psi$ that is differentially sector bounded with sector bound $[0, \Gamma]$ and weight $\Theta$.
\end{enumerate}
Moreover, suppose $n_y=n_\Psi$, $\Gamma \in \s_+^{n_y}$, and $\Theta = \Gamma^{-1}$.
Then, under Assumption \ref{ass:symmetric}, the above four statements are equivalent to the following statement as well.
\begin{enumerate}
\setItemnumber{5}
\item There exists $P \in \s_+^{n_x}$ such that the CTLS \eqref{eq:CTLS} is strongly infinitesimally contracting with norm $\|\cdot\|_P$ and rate $\eta$ for each nonlinearity $\Psi$ that is monotone with upper bound $\Gamma$.
\end{enumerate}
\end{theorem}
\smallskip

The proof of Theorem \ref{thm:lmi_isbmonotone} is shown in Appendix \ref{appendix:proof:lmi}.

\begin{theorem}[Discrete-time sector-bounded Lur'e models]
\label{thm:lmi_isbmonotone_discrete}
Consider the DTLS \eqref{eq:DTLS} with Assumption \ref{ass:C}.
Let $A_\mathrm{cl}$ as in \eqref{def_Acl}, $B_\mathrm{cl}$ as in \eqref{def_Bcl}, $\Theta \in \s_+^{n_\Psi}$, $\Gamma \in \R^{n_y \times n_\Psi}$, $G \coloneqq C^\top \Gamma^\top \Theta$, and $\eta$ be a constant satisfying $0<\eta<1$.
Let $G \coloneqq C^\top \Gamma^\top \Theta$.
Then, the following three statements are equivalent:
\begin{enumerate}
\item There exists $P \in \s_+^{n_x}$ that satisfies
\begin{align}
    \begin{multlined}
    \begin{bmatrix}
        A_\mathrm{cl}^\top P A_\mathrm{cl} - \eta^2 P & A_\mathrm{cl}^\top P B_\mathrm{cl} + G \\
        B_\mathrm{cl}^\top P A_\mathrm{cl} + G^\top & B_\mathrm{cl}^\top P B_\mathrm{cl} - 2 \Theta
    \end{bmatrix}
    \preceq 0.
    \end{multlined}
    \label{contraction_multiplier_discrete}
\end{align}
\item There exists $P \in \s_+^{n_x}$ that satisfies
\begin{align}
    \begin{bmatrix}
        - \eta^2 W & 
        W G &
        (AW + BZ)^\top
        \\
        G^\top W &
        -2 \Theta &
        B_\mathrm{cl}^\top
        \\
        AW + BZ &
        B_\mathrm{cl} &
        - W
    \end{bmatrix}
    \preceq 0
    ,
    \label{contraction_lmi_discrete}
\end{align}
where $W \coloneqq P^{-1}, Z \coloneqq K P^{-1}, B_\mathrm{cl} = B_\Psi + BK_\Psi$.
\item There exists $P \in \s_+^{n_x}$ such that the DTLS \eqref{eq:DTLS} is strongly contracting with norm $\|\cdot\|_P$ and factor $\eta$ for each nonlinearity $\Psi$ that is incrementally sector bounded with sector bound $[0, \Gamma]$ and weight $\Theta$.
\end{enumerate}
Furthermore, if $\Psi$ is differentiable, the above two statements are equivalent to the following statement as well.
\begin{enumerate}
\setItemnumber{4}
\item There exists $P \in \s_+^{n_x}$ such that the DTLS \eqref{eq:DTLS} is strongly contracting with norm $\|\cdot\|_P$ and factor $\eta$ for each nonlinearity $\Psi$ that is differentially sector bounded with sector bound $[0, \Gamma]$ and weight $\Theta$.
\end{enumerate}
Moreover, suppose $n_y=n_\Psi$, $\Gamma \in \s_+^{n_y}$, and $\Theta = \Gamma^{-1}$.
Then, under Assumption \ref{ass:symmetric}, the above three statements are equivalent to the following statement as well.
\begin{enumerate}
\setItemnumber{5}
\item There exists $P \in \s_+^{n_x}$ such that the DTLS \eqref{eq:DTLS} is strongly contracting with norm $\|\cdot\|_P$ and factor $\eta$ for each nonlinearity $\Psi$ that is monotone with upper bound $\Gamma$.
\end{enumerate}
\end{theorem}
\smallskip

The proof of Theorem \ref{thm:lmi_isbmonotone_discrete} is shown in Appendix \ref{appendix:proof:lmi_discrete}.



\subsection{Comparison with LMIs in other literature}
\label{subsection:comparison}

\subsubsection{Lipschitz case}
\label{subsubsection:comparison_lipschitz}

\begin{lemma}
\label{prop:lmi_lipschitz_smaller}
Suppose $n_x = n_y = n_\Psi, B_\Psi = C = \Theta_y = \Theta_\Psi = I_{n_x}, K_\Psi=0$, $\Theta_y=\Theta_\Psi=I_{n_x}$.
Then, the following matrix inequality is sufficient for \eqref{contraction_lmi_lipschitz}:
\begin{align}
    \langle AW + BZ \rangle
    + 2(\eta+\rho) W
    \preceq 0
    .
    \label{contraction_lmi_lipschitz_smaller}
\end{align}
\end{lemma}
\smallskip
\begin{proof}
Substitution of $B_\Psi = C = \Theta_y = \Theta_\Psi = I_{n_x}$ and $K_\Psi=0$ into \eqref{contraction_lmi_lipschitz} leads to the following inequality:
\begin{align}
    \begin{bmatrix}
        \langle AW + BZ \rangle + 2\eta W & 
        I_{n_x} &
        W
        \\
        I_{n_x} &
        - I_{n_x} &
        0
        \\
        W &
        0 &
        - \frac{1}{\rho^2} I_{n_x}
    \end{bmatrix}
    \preceq 0
    ,
    \label{contraction_lmi_lipschitz_substituted}
\end{align}
By taking the Schur complement of the lower right 2$\times$2 block in the left hand side of \eqref{contraction_lmi_lipschitz_substituted}, we observe that \eqref{contraction_lmi_lipschitz_substituted} is equivalent to the following matrix inequality:
\begin{align*}
    \langle AW + BZ \rangle
    + 2\eta W + I_{n_x} + \rho^2 W W
    \preceq 0
    .
\end{align*}
Young's inequality implies $2\rho W
\preceq I_{n_x} + \rho^2 W C^\top C W$.
Therefore, \eqref{contraction_lmi_lipschitz_substituted} implies \eqref{contraction_lmi_lipschitz_smaller}.
\end{proof}

\cite[Lemma 1]{SY-CM-HC-FA:13} claims that robust invariance against an upper-bounded disturbance $w$ is guaranteed by the exponential stability of the error system (i.e., strong infinitesimal contractivity), and an LMI for the controller design that guarantees the contractivity of the closed-loop system is provided in the upper left part of (14) in \cite{SY-CM-HC-FA:13}.
Notably, Lemma \ref{prop:lmi_lipschitz_smaller} tells us that the LMI in \cite{SY-CM-HC-FA:13} is a special case of the LMI \eqref{contraction_lmi_lipschitz} in Theorem \ref{thm:lmi_lipschitz}.
In fact, the negative semidefiniteness of the upper left block of (14) in \cite{SY-CM-HC-FA:13} is identical to \eqref{contraction_lmi_lipschitz_smaller} if the Lipschitz constant $L$ (see (11) in \cite{SY-CM-HC-FA:13}) is renamed as $\rho$ and contraction rate (see (6) in \cite{SY-CM-HC-FA:13}) is redefined as $2\eta$ and if we select a constant $\lambda_0$ as $2(\eta + \rho)$.

\subsubsection{Incremental sector bounded and monotone case}

Propositions 1 and 2 in \cite{MG-VA-ST-DA:23} have shown that \eqref{contraction_multiplier} is sufficient for the contractivity (under a slightly different definition) under the assumption that the Lur'e nonlinearity is differentially sector bounded or monotone, respectively.
Theorem \ref{thm:lmi_isbmonotone} also claims the converse.
In addition, Lemma \ref{lemma:equiv_isb_monotonic} claims that monotonicity is a special case of differential sector boundedness.



\section{EXAMPLE OF CONTROLLER DESIGN}

Consider the DTLS \eqref{eq:DTLS} with the following $A, B, B_\Psi, C$:
\begin{gather}
A =
\begin{bmatrix}
1.2 & 0 & 0 \\
0.1 & 0.8 & 0 \\
0 & 0.1 & 0.6
\end{bmatrix}
,
B =
\begin{bmatrix} 0.2 \\ 0 \\ 0 \end{bmatrix}
,
B_\Psi =
\begin{bmatrix} 0 \\ 0 \\ 0.2 \end{bmatrix}
,
\\
C =
\begin{bmatrix} 1 & 0 & 0 \\ 0 & 1 & 0 \end{bmatrix}
.
\end{gather}
We design gain matrices $K = \begin{bmatrix} k_1 & k_2 & k_3 \end{bmatrix}, K_\Psi \in \R$.
Suppose that the Lur'e nonlinearity $\Psi$ is $\rho$-Lipschitz with $\rho = 0.5$, $\Theta_y = \mathrm{diag} [4, 1]$, $\Theta_\Psi = 1$.
Note that $B_\Psi + B K_\Psi \neq 0$ for all $K_\Psi \in \R$.

Let $\eta = 0.9$.
Then, $W = \mathrm{diag}[0.1, 0.05, 0.2], \ Z = \begin{bmatrix} -0.6 & -0.03 & 0.3 \end{bmatrix}, \ K_\Psi = -1$ is a feasible solution of LMI \eqref{contraction_lmi_lipschitz_discrete} in Theorem \ref{thm:lmi_lipschitz_discrete}, which leads to $K = \begin{bmatrix} -6 & -0.6 & 1.5 \end{bmatrix}$.
Thus, we obtain the controller \eqref{def_controller_discrete} which guarantees the closed-loop contractivity.

To show that the closed-loop system is strongly contracting with these gains, we select three initial states $x^{(1)}(0) = [1;1;1], \ x^{(2)}(0) = [-1;-1;-1]$.
As examples of Lipschitz nonlinearity, we select $\Psi_{(1)}(x^{(1)}, x^{(2)}) = \frac{1}{10} \log(\exp(5x_2)+\exp(-5x_2))+7$, $\Psi_{(2)}(x_1, x_2) = 0.5 / (1+\exp(0.5x_1-x_2))-5$, $\Psi_{(3)}(x_1, x_2) = 0.5 \cos(0.5x_1) \sin(x_2)$.
For each Lur'e nonlinearity and initial condition, we calculate the trajectories $x(k)$ for $k=1,\ldots,10$ according to \eqref{eq:DTLS}.

\begin{figure}[t]
\centering
\includegraphics[width=0.9\linewidth]{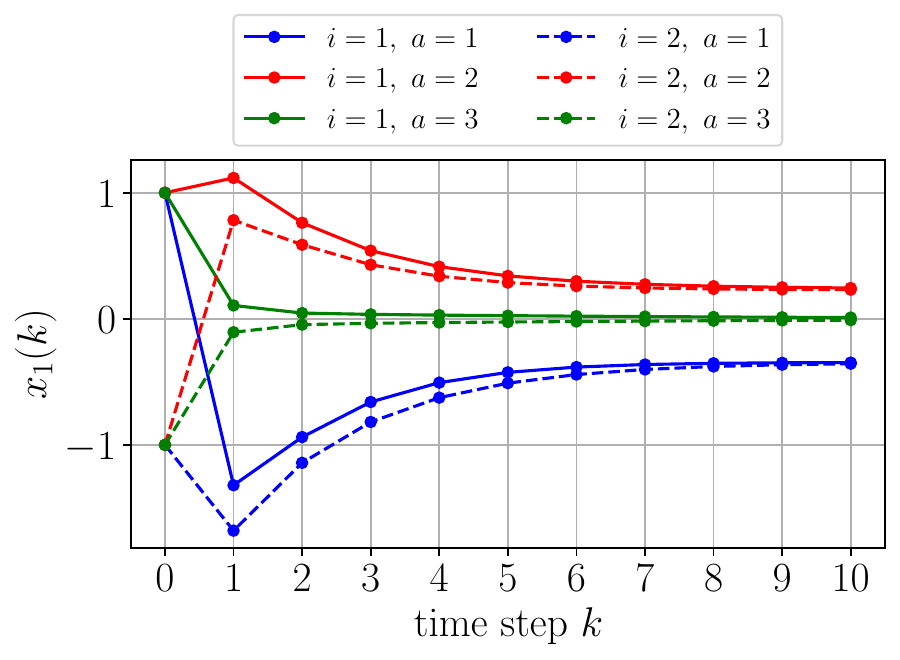}
\caption{Trajectories of $x_1$ for each nonlinearity and initial condition.  \label{fig:trajectories}}
\end{figure}

Each trajectory of the first state $x_1(k)$ for each initial condition $x^{(i)}(0), i=1,2$ and Lur'e nonlinearity $\Psi_{(a)}, a=1,2,3$ is depicted in Fig.~\ref{fig:trajectories}, where the solid lines represent the trajectories for the initial condition $x^{(1)}(0)$ and the dashed lines for $x^{(2)}(0)$, while the blue lines represent the trajectories for nonlinearity $\Psi_{(1)}$, the red lines for $\Psi_{(2)}$, and the green lines for $\Psi_{(3)}$.
It shows that trajectories of $x_1$ with initial states $x^{(1)}(0), x^{(2)}(0)$ converge to the same value for each nonlinearity.
In addition, we evaluate the distance $\| x^{(1)}(k) - x^{(2)}(k) \|_P$, where $P = W^{-1} = \mathrm{diag}[10, 20, 5]$.
The maximum value of the decreasing rate $\| x^{(1)}(k+1) - x^{(2)}(k+1) \|_P/\| x^{(1)}(k) - x^{(2)}(k) \|_P, k=0,\ldots,9$ is observed to be $0.658$, which is less than $\eta=0.9$.

\section{CONCLUSION}

In this paper, we have proposed LMIs for designing controllers that
guarantee contractivity of closed-loop Lur’e systems under Lipschitz,
incrementally sector-bounded, and monotone nonlinearities. We demonstrated
our approach through a practical example illustrating controller synthesis
using the developed LMIs. Our results expand upon previous frameworks,
offering both necessary and sufficient conditions, and enabling broader
applicability to various classes of nonlinear systems. Potential
applications of our contraction-based controller design include tube-based
model predictive control (tube-MPC) \cite{SY-CM-HC-FA:13} and
proportional-integral (PI) control \cite{MG-VA-ST-DA:23}. Future research
directions include further generalizations of incremental quadratic
constraints to a wider range of nonlinear dynamics.

\section*{APPENDIX}
\setcounter{section}{0}
\def\thesection{\Alph{section}}

\section{PROOF OF THEOREM \ref{thm:lmi_lipschitz}}
\label{appendix:proof:lmi_lipschitz}

\begin{proof}
Schur complement of $-\frac{1}{\rho^2} \Theta_y^{-1}$ in the left hand side of \eqref{contraction_lmi_lipschitz} leads to the equivalence between \eqref{contraction_lmi_lipschitz} and the following inequality:
\begin{align}
    \begin{aligned}
    &
    \begin{bmatrix}
        \langle AW + BZ \rangle + 2\eta W + \rho^2 W C^\top \Theta_y C W & 
        B_\mathrm{cl}
        \\
        B_\mathrm{cl}^\top &
        - \Theta_\Psi
    \end{bmatrix}
    \preceq 0
    .
    \end{aligned}
    \notag
\end{align}
Upon pre- and post-multiplying this inequality by $\mathrm{diag}[P, I]$ (which is a congruence transformation and therefore does not lose equivalence), we obtain the equivalence 1) $\Leftrightarrow$ 2).

To prove the equivalence 1) $\Leftrightarrow$ 3), note that
\eqref{contraction_incremental} can be reformulated as follows:
\begin{align*}
    \begin{bmatrix}
        \Delta x \\
        \Delta_x \Psi
    \end{bmatrix}
    ^\top
    \begin{bmatrix}
        \langle P A_\mathrm{cl} \rangle + 2\eta P & P B_\mathrm{cl} \\
        B_\mathrm{cl}^\top P & 0
    \end{bmatrix}
    \begin{bmatrix}
        \Delta x \\
        \Delta_x \Psi
    \end{bmatrix}
    \leq 0
    ,
\end{align*}
where $\Delta x$ and $\Delta_x \Psi$ is defined in \eqref{def_deltax_deltapsi}.
On the other hand, Assumption \ref{ass:C} implies that, for each $y^{(1)}, y^{(2)} \in \R^{n_y}$, there exist $x^{(1)}, x^{(2)} \in \R^{n_x}$ satisfying $y^{(1)}=Cx^{(1)}$, $y^{(2)}=Cx^{(2)}$.
Therefore, \eqref{Lipschitz} holds for any $y^{(1)}, y^{(2)} \in \R^{n_y}$ if and only if the following inequality holds for any $x^{(1)}, x^{(2)} \in \R^{n_x}$:
\begin{align}
    \begin{bmatrix}
        \Delta x \\
        \Delta_x \Psi
    \end{bmatrix}
    ^\top
    \begin{bmatrix}
        \rho^2 C^\top \Theta_y C & 0
        \\
        0 & - \Theta_\Psi
    \end{bmatrix}
    \begin{bmatrix}
        \Delta x \\
        \Delta_x \Psi
    \end{bmatrix}
    \geq 0
    ,
    \label{lipschitz_quadratic}
\end{align}
where $\Delta x$ and $\Delta_x \Psi$ is defined in \eqref{def_deltax_deltapsi}.
Furthermore, the arbitrariness of $x^{(1)}$, $x^{(2)}$, and $\Psi$ implies arbitrariness of $[\Delta x ; \Delta_x \Psi] \in \R^{n_x+n_\Psi}$.
S-lemma states that the claim that ``any vector $[\Delta x ; \Delta_x \Psi] \in \R^{n_x+n_\Psi}$ that satisfies \eqref{Lipschitz} satisfies \eqref{contraction_incremental}'' is equivalent to the existence of $\lambda \geq 0$ such that the following inequality holds:
\begin{align}
    \begin{bmatrix}
        \langle P A_\mathrm{cl} \rangle + 2\eta P & P B_\mathrm{cl} \\
        B_\mathrm{cl}^\top P & 0
    \end{bmatrix}
    +
    \lambda
    \begin{bmatrix}
        \rho^2 \Theta_y & 0
        \\
        0 & - \Theta_\Psi
    \end{bmatrix}
    \preceq 0.
    \label{contraction_multiplier_lipschitz_lambda}
\end{align}
If $\lambda=0$, then \eqref{contraction_multiplier_lipschitz_lambda} implies $PB_\mathrm{cl} = 0$, which contradicts with $B_\mathrm{cl} \neq 0$ and $P \succ 0$.
Therefore, we have $\lambda > 0$.
By redefining $P/\lambda$ as $P$, we obtain \eqref{contraction_multiplier_lipschitz}.
\end{proof}

\section{PROOF OF THEOREM \ref{thm:lmi_lipschitz_discrete}}
\label{appendix:proof:lmi_lipschitz_discrete}

\begin{proof}
Schur complement of the lower right 2$\times$2 block in the left hand side of \eqref{contraction_lmi_lipschitz_discrete} leads to the equivalence between \eqref{contraction_lmi_lipschitz_discrete} and the following inequality:
\begin{align}
    \begin{aligned}
    &
    \begin{bmatrix}
        - \eta^2 W + \rho^2 W C^\top \Theta_y C W & 
        0
        \\
        0 &
        - \Theta_\Psi
    \end{bmatrix}
    \\
    &+
    \begin{bmatrix}
        (AW + BZ)^\top
        \\
        (B_\Psi + BK_\Psi)^\top
    \end{bmatrix}
    P
    \begin{bmatrix}
        (AW + BZ)^\top
        \\
        (B_\Psi + BK_\Psi)^\top
    \end{bmatrix}^\top
    \preceq 0
    ,
    \end{aligned}
    \label{contraction_lmi_lipschitz_discrete_schur}
\end{align}
By pre- and post-multiplying $\mathrm{diag}[P, I]$ by \eqref{contraction_lmi_lipschitz_discrete_schur} (which is a congruence transformation and therefore does not lose equivalence), we obtain the equivalence 1) $\Leftrightarrow$ 2).

We remark that \eqref{contraction_incremental_discrete} can be reformulated as follows:
\begin{align*}
    \begin{bmatrix}
        \Delta x \\
        \Delta_x \Psi
    \end{bmatrix}
    ^\top
    \begin{bmatrix}
        A_\mathrm{cl}^\top P A_\mathrm{cl} - \eta^2 P &
        A_\mathrm{cl}^\top P B_\mathrm{cl} \\
        B_\mathrm{cl}^\top P A_\mathrm{cl} &
        B_\mathrm{cl}^\top P B_\mathrm{cl}
    \end{bmatrix}
    \begin{bmatrix}
        \Delta x \\
        \Delta_x \Psi
    \end{bmatrix}
    \leq 0
    ,
\end{align*}
where $\Delta x$ and $\Delta_x \Psi$ is defined in \eqref{def_deltax_deltapsi}.
Identically to the proof of Theorem \ref{thm:lmi_lipschitz}, S-lemma tells us that 3) is equivalent to the existence of $\lambda \geq 0$ such that the following inequality holds.
\begin{align}
    \begin{bmatrix}
        A_\mathrm{cl}^\top P A_\mathrm{cl} - \eta^2 P &
        A_\mathrm{cl}^\top P B_\mathrm{cl} \\
        B_\mathrm{cl}^\top P A_\mathrm{cl} &
        B_\mathrm{cl}^\top P B_\mathrm{cl}
    \end{bmatrix}
    +
    \lambda
    \begin{bmatrix}
        \rho^2 \Theta_y & 0
        \\
        0 & - \Theta_\Psi
    \end{bmatrix}
    \preceq 0
    .
    \label{contraction_multiplier_discrete_lambda}
\end{align}
If $\lambda=0$, then \eqref{contraction_multiplier_discrete_lambda} implies $B_\mathrm{cl}^\top P B_\mathrm{cl} \preceq 0$, which contradicts with $B_\mathrm{cl} \neq 0$ and $P \succ 0$.
Therefore, we have $\lambda>0$.
Thus, we obtain \eqref{contraction_multiplier_lipschitz_discrete} by redefining $P/\lambda$ as $P$.
\end{proof}

\section{PROOF OF LEMMA \ref{lemma:equiv_isb_monotonic}}
\label{appendix:proof:equiv_isb_monotone}

The statement of Lemma \ref{lemma:equiv_isb_monotonic} is directly deduced by the following lemma with $S = \frac{1}{2} \langle \frac{\partial \Psi}{\partial y}(y) \rangle$.

\begin{lemma}
    \label{lemma:app:equiv}
    Let $\Gamma \in \s_+^n$ and $S \in \R^{n\times n}$ be symmetric.
    Then, the following equivalence holds:
    \begin{align}
        0 \preceq S \preceq \Gamma
        \ \ \Longleftrightarrow \ \ 
        S \Gamma^{-1} (S- \Gamma) \preceq 0
        .
    \label{quadratic_inequality}
    \end{align}
\end{lemma}
\smallskip

\begin{proof}
$S \Gamma^{-1} (S- \Gamma)$ is the Schur's complement of $\Gamma$ in the following matrix:
\begin{align}
\hat{S} \coloneqq
\begin{bmatrix} S & S \\ S & \Gamma \end{bmatrix}
    .
\end{align}
Therefore, what we need to prove is the following equivalence:
\begin{align}
    0 \preceq S \preceq \Gamma
    \ \ \Longleftrightarrow \ \ 
    \hat{S} \succeq 0
    .
\end{align}

($\Leftarrow$) Owing to $S \succeq 0$, we obtain 
\begin{align}
    0 \preceq
    \begin{bmatrix} I \\ I \end{bmatrix}
    S
    \begin{bmatrix} I \\ I \end{bmatrix} ^\top
    =
    \begin{bmatrix} S & S \\ S & S \end{bmatrix}
    .
\end{align}
In addition, $S \preceq \Gamma$ implies the following relation:
\begin{align}
    \begin{bmatrix} S & S \\ S & S \end{bmatrix}
    \preceq
    \begin{bmatrix} S & S \\ S & \Gamma \end{bmatrix}
    = \hat{S}
    .
\end{align}
Therefore, we have $\hat{S} \succeq 0$.

($\Rightarrow$) $S \succeq 0$ comes from $\hat{S} \succeq 0$ because $S$ is the upper left block of $\hat{S}$.
In addition, $\hat{S}\succeq 0$ implies the following relation:
\begin{align}
    0 \preceq
    \begin{bmatrix} I \\ -I \end{bmatrix} ^\top
    \hat{S}
    \begin{bmatrix} I \\ -I \end{bmatrix}
    = \Gamma - S
    ,
\end{align}
which leads us to $S \preceq \Gamma$.
\end{proof}

\section{PROOF OF THEOREM \ref{thm:lmi_isbmonotone}}
\label{appendix:proof:lmi}

\begin{proof}
On pre- and post-multiplying $\mathrm{diag}[P, I]$ by \eqref{contraction_lmi} (which is a congruence transformation and therefore does not lose equivalence), we obtain the equivalence 1) $\Leftrightarrow$ 2).

Regarding the equivalences 1) $\Leftrightarrow$ 3), we observe that, providing Assumption \ref{ass:C}, \eqref{incremental_sector_boundedness} holds for any $y^{(1)}, y^{(2)} \in \R^{n_y}$ if and only if the following inequality holds for any $x^{(1)}, x^{(2)} \in \R^{n_x}$:
\begin{align}
    \begin{bmatrix}
        \Delta x \\
        \Delta_x \Psi
    \end{bmatrix}
    ^\top
    \begin{bmatrix}
        0 & C^\top \Gamma^\top \Theta \\
        \Theta \Gamma C & - 2 \Theta
    \end{bmatrix}
    \begin{bmatrix}
        \Delta x \\
        \Delta_x \Psi
    \end{bmatrix}
    \geq 0
    ,
    \label{isb_quadratic}
\end{align}
where $\Delta x$ and $\Delta_x \Psi$ is defined in \eqref{def_deltax_deltapsi}.
The rest of the proof is parallel to the proof of Theorem \ref{thm:lmi_lipschitz}.

To prove the equivalence 1) $\Leftrightarrow$ 4), we consider an arbitrary vector $\delta x \in \R^{n_x}$ and introduce the following vector:
\begin{align}
    \delta v
    & \coloneqq 
    \begin{bmatrix}
        I \\\frac{\partial \Psi}{\partial y}(Cx) C
    \end{bmatrix}
    \delta x
    .
\end{align}
We remark that \eqref{contraction_differential} holds if and only if the following inequality holds for any $\delta x \in \R^{n_x}$:
\begin{align}
    \delta v^\top
    &
    \begin{bmatrix}
        \langle P A_\mathrm{cl} \rangle + 2\eta P & P B_\mathrm{cl} \\
        B_\mathrm{cl}^\top P & 0
    \end{bmatrix}
    \delta v
    \leq 0
    .
    \label{statement_differential}
\end{align}
On the other hand, Assumption \ref{ass:C} implies that there exists $x \in \R^{n_x}$ such that $y = Cx$ and that $C \delta x$ spans $\R^{n_y}$.
Therefore, on pre- and post-multiplying $C \delta x$ by \eqref{differentially_sector_bounded} and reformulating the product, we observe that \eqref{differentially_sector_bounded} is equivalent to the following statement:
\begin{align}
    \delta v^\top
    \begin{bmatrix}
        0 & C^\top \Gamma^\top \Theta \\
        \Theta \Gamma C & - 2 \Theta
    \end{bmatrix}
    \delta v
    \geq 0
    \ \ \ 
    \forall \delta x\in \R^{n_x}
    .
\end{align}
Furthermore, the arbitrariness of $\delta x$ and $\Psi$ implies the arbitrariness of $\delta v$.
Therefore, S-lemma implies that ``\eqref{differentially_sector_bounded} implies \eqref{contraction_differential}'' is equivalent to the existence of $\lambda \geq 0$ that satisfies the following inequality:
\begin{align}
    \begin{bmatrix}
        \langle P A_\mathrm{cl} \rangle + 2\eta P & P B_\mathrm{cl} \\
        B_\mathrm{cl}^\top P & 0
    \end{bmatrix}
    +
    \lambda
    \begin{bmatrix}
        0 & C^\top \Gamma^\top \Theta \\
        \Theta \Gamma C & - 2 \Theta
    \end{bmatrix}
    \preceq 0.
    \notag
\end{align}
As in the proof of Theorem \ref{thm:lmi_lipschitz}, $B_\mathrm{cl} \neq 0$ leads to $\lambda > 0$, and again we obtain \eqref{contraction_multiplier} by redefining $P/\lambda$ as $P$.

Lemma \ref{lemma:equiv_isb_monotonic} proves the equivalence 4) $\Leftrightarrow$ 5).
\end{proof}

\section{PROOF OF THEOREM \ref{thm:lmi_isbmonotone_discrete}}
\label{appendix:proof:lmi_discrete}

\begin{proof}
Schur complement of $-W$ in the left hand side of \eqref{contraction_lmi_discrete} leads to the equivalence between \eqref{contraction_lmi_discrete} and the following inequality:
\begin{align}
    \begin{aligned}
    &
    \begin{bmatrix}
        - \eta^2 W & 
        WC^\top \Gamma^\top \Theta
        \\
        \Theta \Gamma C W &
        -2 \Theta
    \end{bmatrix}
    \\
    &+
    \begin{bmatrix}
        (AW + BZ)^\top
        \\
        (B_\Psi + BK_\Psi)^\top
    \end{bmatrix}
    P
    \begin{bmatrix}
        (AW + BZ)^\top
        \\
        (B_\Psi + BK_\Psi)^\top
    \end{bmatrix}^\top
    \preceq 0
    ,
    \end{aligned}
    \label{contraction_lmi_discrete_schur}
\end{align}
On pre- and post-multiplying $\mathrm{diag}[P, I]$ by \eqref{contraction_lmi_discrete_schur} (which is a congruence transformation and therefore does not lose equivalence), we obtain the equivalence 1) $\Leftrightarrow$ 2).

The proof of the equivalence 1) $\Leftrightarrow$ 3) is parallel to that in the proof of Theorem \ref{thm:lmi_lipschitz_discrete} except that Lipschitzness is replaced by incremental sector boundedness.

The proof of the equivalence 1) $\Leftrightarrow$ 4) is parallel to that in the proof of Theorem \ref{thm:lmi_isbmonotone} except that the contractivity of the CTLS is replaced by that of the DTLS.

Lemma \ref{lemma:equiv_isb_monotonic} proves the equivalence 4) $\Leftrightarrow$ 5).
\end{proof}

\bibliographystyle{plainurl+isbn}
\bibliography{alias,Main,FB}





\end{document}